\def\hldest#1#2#3{}
% Macros for creating documents
% Written by Marcel Goh, unless otherwise specified.

% OLD FONT NAMES. These are not used anymore (as of 29 Jul 2020),
% but left here for compatibility purposes.

\font\smallheader=cmssbx10 % for section headers

% Roman fonts
\font\eightpt=cmr8
\font\ninept=cmr9

% Bold fonts

\font\mathbold=cmmib10

% NEW FONT NAMES FOR PRELOADED PLAIN FONTS (naming mostly consistent with the TeXbook)
\font\ninerm=cmr9     \font\eightrm=cmr8   \font\sixrm=cmr6      % ROMAN
\font\ninei=cmmi9     \font\eighti=cmmi8   \font\sixi=cmmi6      % MATH ITALIC
\font\ninesy=cmsy9    \font\eightsy=cmsy8  \font\sixsy=cmsy6     % SYMBOLS
\font\ninebf=cmbx9    \font\eightbf=cmbx8  \font\sixbf=cmbx6     % BOLD EXTENDED
\font\ninett=cmtt9    \font\eighttt=cmtt8                        % TYPEWRITER
\font\nineit=cmti9    \font\eightit=cmti8     % TEXT ITALIC
\font\ninesl=cmsl9    \font\eightsl=cmsl8                        % SLANTED

\font\tensc=cmcsc10   \font\ninesc=cmcsc9  \font\eightsc=cmcsc8  % SMALL CAPS
        % SANS SERIF
    % SANS SERIF ITALIC

\font\eightssq=cmssq8  \font\eightssqi=cmssqi8  % SANS SERIFS FOR QUOTES
 % SLANTED TYPEWRITER
       % UNSLANTED TEXT ITALIC

    % BOLD MATH ITALIC
   % BOLD MATH SYMBOLS

\font\tenssbx=cmssbx10 % SANS SERIF BOLD EXTENDED

% MISC. OTHER FONTS
  \font\twelvebf=cmbx12
  
\def\sc{\tensc}  \def\mc{\ninerm}

% CYRILLIC
\input cyracc.def
    \font\tencyr=wncyr10   \font\ninecyr=wncyr9   \font\eightcyr=wncyr8
    \font\tencyri=wncyi10  \font\ninecyri=wncyi9  \font\eightcyri=wncyi8
    \def\cyr{\tencyr\cyracc} \def\cyri{\tencyri\cyracc}

% FONT SIZES (ADAPTED FROM THE TEXBOOK).
% CHANGES MADE:
% - Edited \sc and added \mc.  (1)
% - Cyrillic support built in. (2)
\newskip\ttglue
\def\tenpoint{\def\rm{\fam0\tenrm}%
  \textfont0=\tenrm \scriptfont0=\sevenrm \scriptscriptfont0=\fiverm
  \textfont1=\teni  \scriptfont1=\seveni  \scriptscriptfont1=\fivei
  \textfont2=\tensy \scriptfont2=\sevensy \scriptscriptfont2=\fivesy
  \textfont3=\tenex \scriptfont3=\tenex   \scriptscriptfont3=\tenex
  \textfont\itfam=\tenit  \def\it{\fam\itfam\tenit}%
  \textfont\slfam=\tensl  \def\sl{\fam\slfam\tensl}%
  \textfont\ttfam=\tentt  \def\tt{\fam\ttfam\tentt}%
  \textfont\bffam=\tenbf  \scriptfont\bffam=\sevenbf
   \scriptscriptfont\bffam=\fivebf \def\bf{\fam\bffam\tenbf}%
  \tt \ttglue=.5em plus.25em minus.15em
  \normalbaselineskip=12pt
  \setbox\strutbox=\hbox{\vrule height8.5pt depth3.5pt width0pt}%
  \let\sc=\tensc \let\mc=\ninerm  % (1)
  \def\cyr{\tencyr\cyracc}\def\cyri{\tencyri\cyracc}% (2)
  \let\big=\tenbig  \normalbaselines\rm}

\def\ninepoint{\def\rm{\fam0\ninerm}%
\textfont0=\ninerm  \scriptfont0=\sixrm  \scriptscriptfont0=\fiverm
\textfont1=\ninei   \scriptfont1=\sixi   \scriptscriptfont1=\fivei
\textfont2=\ninesy  \scriptfont2=\sixsy  \scriptscriptfont2=\fivesy
\textfont3=\tenex   \scriptfont3=\tenex  \scriptscriptfont3=\tenex
\textfont\itfam=\nineit  \def\it{\fam\itfam\nineit}%
\textfont\slfam=\ninesl  \def\sl{\fam\slfam\ninesl}%
\textfont\ttfam=\ninett  \def\tt{\fam\ttfam\ninett}%
\textfont\bffam=\ninebf  \scriptfont\bffam=\sixbf
\scriptscriptfont\bffam=\fivebf\def\bf{\fam\bffam\ninebf}%
\tt\ttglue=.5em plus.25em minus.15em
\normalbaselineskip=11pt
\setbox\strutbox=\hbox{\vrule height8pt depth3pt width0pt}%
\let\sc=\ninesc\let\mc=\eightrm%(1)
\def\cyr{\ninecyr\cyracc}\def\cyri{\ninecyri\cyracc}%(2)
\let\big=\ninebig\normalbaselines\rm}

\def\eightpoint{\def\rm{\fam0\eightrm}%
  \textfont0=\eightrm \scriptfont0=\sixrm \scriptscriptfont0=\fiverm
  \textfont1=\eighti  \scriptfont1=\sixi  \scriptscriptfont1=\fivei
  \textfont2=\eightsy \scriptfont2=\sixsy \scriptscriptfont2=\fivesy
  \textfont3=\tenex   \scriptfont3=\tenex \scriptscriptfont3=\tenex
  \textfont\itfam=\eightit  \def\it{\fam\itfam\eightit}%
  \textfont\slfam=\eightsl  \def\sl{\fam\slfam\eightsl}%
  \textfont\ttfam=\eighttt  \def\tt{\fam\ttfam\eighttt}%
  \textfont\bffam=\eightbf  \scriptfont\bffam=\sixbf
  %\scriptscriptfont\bffam=\fivebf  \def\bf{\fam\bffam\eightbf}
  % \tt \ttglue=.5em plus.25em minus.15em
  \normalbaselineskip=9pt
  %\setbox\strutbox=\hbox{\vrule height7pt depth2pt width0pt}
  \let\sc=\eightsc \let\mc=\sevenrm  % (1)
  \def\cyr{\eightcyr\cyracc}\def\cyri{\eightcyri\cyracc}% (2)
  \let\big=\eightbig  \normalbaselines\rm}%
\def\nospace{\nulldelimiterspace0pt\mathsurround0pt}%
\def\tenbig#1{{\hbox{$\left#1\vbox to8.5pt{}\right.\nospace$}}}%
\def\ninebig#1{{\hbox{$\textfont0=\tenrm\textfont2=\tensy
  \left#1\vbox to7.25pt{}\right.\nospace$}}}%
\def\eightbig#1{{\hbox{$\textfont0=\ninerm\textfont2=\ninesy
  \left#1\vbox to6.5pt{}\right.\nospace$}}}%

% USE NON-EXTENDED VERSIONS OF BOLD FONT
\def\nonextendedbold{
  \font\fiveb=cmb10 at 5pt
  \font\sixb=cmb10 at 6pt
  \font\sevenb=cmb10 at 7pt
  \font\eightb=cmb10 at 8pt
  \font\nineb=cmb10 at 9pt
  \font\tenb=cmb10
  \font\twelveb=cmb10 at 12pt
  \let\fivebf=\fiveb
  \let\sixbf=\sixb
  \let\sevenbf=\sevenb
  \let\eightbf=\eightb
  \let\ninebf=\nineb
  \let\tenbf=\tenb
  \let\twelvebf=\twelveb
}

% Create header on non-first pages
\def\leftrighttop#1#2{
  \headline{\ifnum\pageno=1\hfil\else{\ninept #1 \hfil #2}\fi}
}

% No page numbers on first page
\def\firstnopagenum{
  \footline{\ifnum\pageno=1 \hfil \else \hfil{\rm \number\pageno}\hfil\fi}
}

% Create title, subtitle, name, and date
\def\maketitle#1#2#3#4{
  \centerline {\titlefont #1}
  \medskip
  \centerline {\eightpt #2}
  \medskip
  \centerline {\tensc #3}
  \medskip
  \centerline {\tensc #4}
  \bigskip
}
% Two-author version

% Create title, subtitle, and name

% Create floating text box with specified width in inches.
% Useful for abstracts (6 true inches is good):
% \floattext 6 Abstract. text...
\outer\def\floattext#1 #2. #3\par{
  $$
  \vbox{
    \hsize #1 true in
    \noindent{\bf #2.}\enskip #3
  }
  $$
}

% Math version of floattext, no bold label

% Section header in heavy sans-serif font
% Separate two of these in a row with \vskip -\medskipamount
\def\lsection#1\par{
  \bigskip\vskip\parskip
  \leftline{\sectionfont#1}\nobreak\medskip\noindent
}

% Centred version
\def\csection#1\par{
  \bigskip\vskip\parskip
  \centerline{\sectionfont#1}\nobreak\medskip\noindent
}

% Right-justified version
\def\rsection#1\par{
  \bigskip\vskip\parskip
  \rightline{\sectionfont#1}\nobreak\medskip\noindent
}
\def\section{\lsection}

% New subsection
\def\boldlabel#1. {\noindent{\bf #1.\enspace}}
\def\subsection#1. {\medskip\noindent{\bf #1.\enspace}}

% Block of left-justified text, in "poetry-mode"

% Block of left-justified text, in "poetry-mode"

% Fraktur font -- needs some extra tuning to work with sizes other than \tenpoint
\font\tenfrak=eufm10
\font\sevenfrak=eufm7
\font\fivefrak=eufm5
\newfam\frakfam
\textfont\frakfam=\tenfrak
\scriptfont\frakfam=\sevenfrak
\scriptscriptfont\frakfam=\fivefrak
\def\frak#1{{\fam\frakfam #1}}

% A janky small-caps
\def\janksc#1#2 {#1{\eightpt#2}}
\def\jankscsp#1#2 {#1{\eightpt#2}\ }
\def\scproclaim#1.#2\par{\noindent\jankscsp #1.\enspace{\it#2\par}}

% References
\def\ref#1{[#1]}

% Quotes -- adapted from manmac
\def\quote{
  \begingroup
    \baselineskip 10pt
    \parfillskip 0pt
    \interlinepenalty 10000
    \leftskip 0pt plus 40pc minus \parindent
    \let\rm=\quoterm\let\sl=\quotesl\everypar{\sl}
    \obeylines
}
% This is used to end a \quote block. (I'm not sure why the second \rm is needed.)
\def\author#1(#2){\nobreak\smallskip\rm--- \rm#1\unskip\enspace(#2)\par\endgroup}

% DEFAULT SETTINGS
\def\titlefont{\twelvebf}
\def\sectionfont{\tenssbx}
\def\quoterm{\eightssq}
\def\quotesl{\eightssqi}

% CLASSIC MODE:
% Section headers are bold-serif; bolds and sans-serifs are more condensed

% WIDE MARGINS
\def\widemargins{
  \magnification=\magstephalf\hoffset=40pt \voffset=28pt
  \hsize=29pc  \vsize=45pc  \maxdepth=2.2pt  \parindent=19pt
}
\def\bookheader#1#2{
  \nopagenumbers
  \def\leftheadline{{\rm\folio}\hfil{\eightpoint#1}\hfil}
  \def\rightheadline{\hfil{\eightpoint#2}\hfil{\rm\folio}}
  \headline{\ifodd\pageno{\ifnum\pageno<2\hfil\else\rightheadline\fi}\else\leftheadline\fi}
}

\tenpoint

% Macros for typesetting math (assignments)
% Written by Marcel Goh, except for the parts that are not.

% ================== KNUTH ================== %

\def\xskip{\hskip 7pt plus 3pt minus 4pt}

\def\proof{\medbreak\noindent{\it Proof.}\xskip\ignorespaces}

\def\slug{\quad\hbox{\kern1.5pt\vrule width2.5pt height6pt depth1.5pt\kern1.5pt}\medskip}
\def\noskipslug{\quad\hbox{\kern1.5pt\vrule width2.5pt height6pt depth1.5pt\kern1.5pt}}

% Algorithms
\newdimen\algindent
\newif\ifitempar \itempartrue % normally true unless briefly set false
\def\algindentset#1{\setbox0\hbox{{\bf #1.\kern.25em}}\algindent=\wd0\relax}
\def\algbegin #1 #2{\algindentset{#21}\alg #1 #2} % when steps all have 1 digit
\def\aalgbegin #1 #2{\algindentset{#211}\alg #1 #2} % when 10 or more steps
\def\alg#1(#2). {\medbreak % Usage: \algbegin Algorithm A (algname). This...
  \noindent{\bf#1}({\it#2\/}).\xskip\ignorespaces}
\def\algstep#1.{\ifitempar\smallskip\noindent\else\itempartrue
  \hskip-\parindent\fi
  \hbox to\algindent{\bf\hfil #1.\kern.25em}%
  \hangindent=\algindent\hangafter=1\ignorespaces}

% ================ END KNUTH ================ %

% Sets of numbers (not compatible with CWEB)
\def\NN{{\bf N}}
\def\ZZ{{\bf Z}}

\def\RR{{\bf R}}

  % Deprecated

\def\C{{\bf C}}

% Probability
\def\pr{\mathop{\hbox{\bf P}}\nolimits}
\def\ex{\mathop{\hbox{\bf E}}\nolimits}
\def\var{\mathop{\hbox{\bf V}}\nolimits}
\def\one{\mathop{\hbox{\bf 1}}\nolimits}
\def\indic#1{\one_{\sevenrm[#1]}}

% Graphs

% Analysis

% Counters
\newcount\thmcount  % Counts theorems/lemmas/corollaries
\thmcount=1
\newcount\sectcount  % Counts sections
\sectcount=1
\newcount\figcount  % Counts figures
\figcount=1
\newcount\eqcount  % Counts equations
\eqcount=1

% Number an equation oldstyle
\def\oldno#1{\eqno({\oldstyle#1})}
\def\refeq#1{({\oldstyle#1})}
\def\adveq{\oldno{\the\eqcount}\global\advance\eqcount by 1}  % Print and advance equation counter
\def\advthm{\the\thmcount\global\advance \thmcount by 1}

% Numbered section (advances equation counter)
\def\advsect{\section\the\sectcount\global\advance\sectcount by 1. }

% Proclaim with parentheses in italics
\outer\def\parenproclaim #1 (#2).#3\par{\medbreak
  \noindent{\bf #1}\enspace\rm({\it #2\/}).\nobreak\ignorespaces{\sl #3\par}
  \ifdim\lastskip<\medskipamount \removelastskip\penalty55\medskip\fi}

% Exactly what it sounds like

\input eplain

\def\divides{\backslash}
\def\zed#1{\ZZ/#1\ZZ}
\def\tot{\varphi}
\def\C{{\cal C}}

\def\ref#1{[\hlstart{name}{}{bib#1}#1\hlend]}%Fix what eplain messed up
\def\noblueboxes{\special{ps:[/pdfm { /big_fat_array exch def big_fat_array 1 get 0
0 put big_fat_array 1 get 1 0 put big_fat_array 1 get 2 0 put big_fat_array pdfmnew } def}}

\enablehyperlinks

% pdf outline? 
\ifpdf% not our case [pdftex, other pdf drivers]
  \hlopts{bwidth=0}% no border around the link 
  \pdfoutline goto name {sec1} count -1 {Introduction}%
    \pdfoutline goto name {defsnot} {Definitions and notation}%
  \pdfoutline goto name {sec2} count -1 {Torsion-free groups}%
    \pdfoutline goto name {countexact} {Counting exact values}%
  \pdfoutline goto name {sec3} count -4 {Finite additive groups}%
    \pdfoutline goto name {cyclic} {Cyclic groups}%
    \pdfoutline goto name {exactval} {Exact values}%
    \pdfoutline goto name {arbfinite} {Arbitrary finite abelian groups}%
    \pdfoutline goto name {elempgroups} {Elementary p-groups}%
  \pdfoutline goto name {sec4} {Further work}%
  \pdfoutline goto name {acks} {Acknowledgements}%
  \pdfoutline goto name {refs} {References}%
\else% our case
  \special{ps:[/PageMode /UseOutlines /DOCVIEW pdfmark}%
  % The individual outline entries, using a different syntax than pdftex, but the same information.
  \special{ps:[/Count -1 /Dest (sec1) cvn /Title (Introduction) /OUT pdfmark}%
    \special{ps:[/Count -0 /Dest (defsnot) cvn /Title (Definitions and notation) /OUT pdfmark}%
  \special{ps:[/Count -1 /Dest (sec2) cvn /Title (Torsion-free groups) /OUT pdfmark}%
    \special{ps:[/Count -0 /Dest (countexact) cvn /Title (Counting exact values) /OUT pdfmark}%
  \special{ps:[/Count -4 /Dest (sec3) cvn /Title (Finite additive groups) /OUT pdfmark}%
    \special{ps:[/Count -0 /Dest (cyclic) cvn /Title (Cyclic groups) /OUT pdfmark}%
    \special{ps:[/Count -0 /Dest (exactval) cvn /Title (Exact values) /OUT pdfmark}%
    \special{ps:[/Count -0 /Dest (arbfinite) cvn /Title (Arbitrary finite abelian groups) /OUT pdfmark}%
    \special{ps:[/Count -0 /Dest (elempgroups) cvn /Title (Elementary p-groups) /OUT pdfmark}%
  \special{ps:[/Count -0 /Dest (sec4) cvn /Title (Further work) /OUT pdfmark}%
  \special{ps:[/Count -0 /Dest (acks) cvn /Title (Acknowledgements) /OUT pdfmark}%
  \special{ps:[/Count -0 /Dest (refs) cvn /Title (References) /OUT pdfmark}%
\fi
\noblueboxes

\bookheader{ARITHMETIC SUBSEQUENCES}{MARCEL K. GOH AND ROSIE Y. ZHAO}
\widemargins
\nopagenumbers

\centerline{\bf ARITHMETIC SUBSEQUENCES IN A RANDOM ORDERING}
\centerline{\bf OF AN ADDITIVE SET}
\bigskip
\bigskip
\centerline{\sc Marcel K. Goh}
\smallskip
\centerline{\sl Department of Mathematics and Statistics, McGill University}
\bigskip
\centerline{\sc Rosie Y. Zhao}
\smallskip
\centerline{\sl School of Computer Science, McGill University}
\bigskip

\floattext5 \ninepoint\bf Abstract. \ninepoint
For a finite set $A$ of size $n$, an ordering is an
injection from $\{1,2,\ldots,n\}$ to $A$.
We present results concerning the asymptotic properties of the length
$L_n$ of the longest
arithmetic subsequence in a random ordering of an additive set $A$. In the
torsion-free case
where $A = [1,n]^d\subseteq {\bf Z}^d$, we prove that $L_n\sim 2d\log n/\log\log n$.
We show that the case $A = \zed n$ behaves asymptotically like the torsion-free case with $d=1$, and then
use this fact to compute the expected length of the longest arithmetic subsequence in
a random ordering of an arbitrary finite abelian group. We also prove that the number of orderings of
$\zed n$ without any arithmetic subsequence of length $3$ is $2^{n-1}$ when $n\geq 2$ is a power of $2$,
and zero otherwise. We conclude with a concrete application to elementary
$p$-groups and a discussion of possible noncommutative generalisations.

\advsect Introduction
\hldest{xyz}{}{sec1}

{\sc Arithmetic progressions} are most often viewed as subsets of an ambient abelian group $Z$.
Under this interpretation and with $Z=\ZZ$, the group of integers,
E.~Szemer\'edi's celebrated theorem states that any dense subset $A$ of $Z$ contains arbitrarily long arithmetic
progressions \ref{13}. The story is different when the elements of $A$ are ordered in a sequence
and we look for progressions in {\it subsequences}. For example,
it is always possible to permute the first $n$ integers so that no subsequence of length 3 is an arithmetic
progression~\ref{7}. Using this fact, J. A. Davis, R. C. Entringer,
R.~L.~Graham, and G.~J.~Simmons showed in 1977 that there are permutations of the positive
integers that do
not contain any arithmetic progressions of length 5 and permutations of $\ZZ$ avoiding 7-term arithmetic
progressions~\ref{5}.
Whether or not there exist permutations of the positive
integers not admitting any permutation of length 4 remains an open problem to this day, but
J.~Geneson has recently improved the result of Davis et al.\ for $\ZZ$
to 6-term arithmetic progressions~\ref{9}.

In an arbitrary abelian group $Z$, an entirely different formulation of the problem studies bijections $\sigma$
from $Z$ to itself, declaring a $k$-term arithmetic
progression to be a tuple $(z_1,z_2,\ldots,z_k)\in Z^k$ such that there exists some $r$ for which
$z_{i+1} - z_i = r$ for every $1\leq i<k$. A permutation is said to {\it destroy}\/
($3$-term) arithmetic progressions
if for every triple $(x,y,z)\in Z^3$, either $(x,y,z)$ is not a progression or
$\big(\sigma(x), \sigma(y), \sigma(z)\big)$
is not a progression. For $Z$ infinite, P.~Hegarty showed that there exists a progression-destroying
permutation of $Z$ if and only if $Z/\Omega_2(Z)$ is equipotent to $Z$, where $\Omega_2(Z)$ is the subgroup
of elements with even order~\ref{10}. When $Z = (\zed p)^d$ for
$p$ prime, N.~D.~Elkies and A.~A.~Swaminathan proved that there exist progression-destroying permutations if and
only if $(p,d)\notin\big\{(3,1),(5,1),(7,1)\big\}$ and $p$ is odd~\ref{6}.

The above definition of a progression in a permutation is not the one that we will study. Instead, we will
{\it order} the elements of an additive set $A$ from $1$ to $|A|$ and count
arithmetic progressions that appear in order in the sequence.
(Thus our definition is closer to the one used by Davis et al.;
Hegarty's definition above is akin to requiring that the indices
of an arithmetic subseqence form a progression as well.)
Though it may be possible to find perverse orderings of an additive set $A$ that admit no arithmetic
subsequences of a given length $k$,
it is natural to suspect that as the size of $A$ gets large, this becomes more and more unlikely in some
formal sense. We can then replace a difficult extremal problem with a simpler probabilistic one.
In this paper, we show that when $A = \{1,2,\ldots,n\}\subseteq \ZZ$,
and the ordering is taken
uniformly at random from $\frak S_n$ (the group of permutations of an $n$-element set),
the length of the longest arithmetic progression
converges to $2\log n/\log\log n$ in probability. We also show analogous results when $Z$ is a finite
cyclic group and $A$ is taken to be the whole group. In the case of an arbitrary finite abelian group $Z$,
one cannot
make any sensible asymptotic statements that depend only on the order of the group, but we will see that
the counting methods used in the $\ZZ^d$ and cyclic cases can be extended to give bounds on the expected number of
arithmetic subsequences of a given length $k$ in a random ordering of $Z$.

\noblueboxes
Results of this type have appeared in the literature
for progressions (view\-ed as subsets) in a random set. In the case
of a random subset of $\{1,2,\ldots,n\}$, where each
element is included independently with probability $1/2$, I.~Benjamini, A.~Yadin, and O.~Zeitouni showed
that the expected length of the longest arithmetic progression is asymptotically $2\log n/\log 2$, and
that the asymptotics do not change if the integers are taken modulo $n$~\ref{2}.
Returning to the setting of random sequences, A.~Frieze proved in 1991 that the length of the
longest monotone increasing subsequence of a random permutation of $\{1,2,\ldots,n\}$ is concentrated around
$2\sqrt n$, its asymptotic expected value~\ref{8}.

% Tangential
%On a somewhat tangential note, finding the longest arithmetic subsequence in a given sequence
%(with values taken from a finite set) is a textbook example of a programming problem
%for which the technique of dynamic programming yields an efficient solution.
%This makes the ``longest arithmetic subsequence'' a common interview question in the
%technology industry.

\medskip
\hldest{xyz}{}{defsnot}
\boldlabel Definitions and notation.
We will make frequent use of the following asymptotic notation for nonnegative-valued
functions $f,g:\NN\to\RR$.
We write $f(n) = O_n\big(g(n)\big)$ if there exists a constant $C$ such that $f(n)\leq C\cdot g(n)$ for all
$n\in \NN$. When $\lim_{n\to\infty} f(n)/g(n) = 0$, we write $f(n)=o_n\big(g(n)\big)$, and if instead we
have $\lim_{n\to\infty} f(n)/g(n) = \infty$, then we write $f(n)=\omega_n\big(g(n)\big)$. Lastly, if
$\lim_{n\to\infty} f(n)/g(n) = 1$, then we write $f(n)\sim_n g(n)$. In all cases, we may omit the subscript
when the asymptotic variable is clear from context.

In this paper, the notation $[a,b]$ will indicate the {\it discrete} interval $\{n\in \ZZ : a\leq n\leq b\}$,
unless otherwise stated.
Formally, an {\it additive set} is a pair $(A,Z)$, where $(Z,+)$ is an
abelian group
with identity $0$
and $A\subseteq Z$ is finite and nonempty. When the ambient group is evident from context,
we will write $A$ for $(A,Z)$.
For $n\in \ZZ$ and $z\in Z$, we let $nz$ denote the iterated sum; so $2z = z+z$ and $-3z = -z-z-z$. The {\it
order} $|z|$ of a group element $z$ is the smallest positive integer $n$ such that $nz = 0$.
If no such integer exists, then we put $|z| = \infty$. A group all of whose elements have finite order is called
a {\it torsion group}, while a group in which every non-identity element has infinite order is a {\it torsion-free
group}.

For our purposes, an {\it ordering} $\sigma$ of an additive set $A$ is a bijection $\sigma : [1,|A|] \to A$
and we will regard such an ordering as a sequence by listing its elements
$\sigma = \big(\sigma(1),\sigma(2),\ldots,\sigma(|A|)\big)$ in order.
%(We will try to avoid using the
%term {\it permutation}, which is formally a bijection from a set $A$ to itself. Of course, when $A$ is
%the set $[1,n]$, then an ordering is indeed a {\it bona fide} permutation, and in all other cases the two
%concepts are largely the same, so the reader may mentally substitute the word ``permutation'' whenever
%the word ``ordering'' appears, without any detriment to comprehension.)
More generally, a {\it $k$-ordering} of a set $A$, where $k\leq|A|$, is an injection from $[1,k]\to |A|$
and it can also be regarded as a sequence.
We will say that a subsequence
$\sigma_{i_1},\sigma_{i_2},\ldots,\sigma_{i_k}$ with $1\leq i_1<i_2<\cdots<i_k\leq n$
is an {\it arithmetic progression} (or simply {\it progression}) of length $k$ if for some $a,r\in Z$,
$$(\sigma_{i_1},\sigma_{i_2},\ldots,\sigma_{i_k}) = \big(a, a+r, \ldots, a+(k-1)r\big).$$
The element $a$ is called the {\it base point} and $r$ is called the {\it step} or {\it common difference}.
If $r=0$, the progression is said to be {\it trivial}\/; because we will most often be looking at progressions
in sequences without repetition, a progression will be assumed to be nontrivial unless otherwise stated.

Let $\sigma$ be an ordering of an additive set $A$ of size $n$. We define
$L(\sigma)$ to be the largest $k$ for which there exists a subsequence of $\sigma$ that is an arithmetic
progression of length $k$. For example, if $\sigma$ is the ordering $(2,7,1,6,3,4,5)$ of $[1,7]$, then
$L(\sigma) = 4$, due to the embedded 4-ordering $(2,3,4,5)$.
We can make $\sigma$ random by fixing an enumeration $f:[1,n]\to A$, selecting
$\pi$ uniformly at random from the group $\frak S_n$, and then letting $\sigma = f\circ \pi$.
This makes $L(\sigma)$ a random variable, which we will denote by $L_n$. We are interested in determining
how $L_n$ grows as $n$ gets large.

\hldest{xyz}{}{sec2}
\advsect Torsion-free groups

By the structure theorem of finitely-generated abelian torsion-free groups, every abelian torsion-free group
is isomorphic to a lattice $\ZZ^d$ for $d\geq 0$ (when $d=0$, this is the trivial group, so we will assume
that $d\geq1$ for the rest of this section).
A natural subset of $\ZZ^d$ to consider is $A = [1,n]^d$.
We first prove a lemma for the simplest case $d=1$, then generalise the result to any value of $d$.

\newcount\kperms
\newcount\kpermseq
\newcount\kpermsineq
\kperms=\thmcount
\proclaim Lemma \advthm. For $1\leq k\leq n$, let
$P_{nk}$ be the number of $k$-orderings of $[1,n]$ that are arithmetic progressions, i.e., of the form
$$\big(a, a+r, a+2r, \ldots, a+(k-1)r\big)$$
for some $a,r\in\ZZ$. We have $P_{n1} = n$ and
\global\kpermseq=\eqcount
$$P_{nk} = 2n\left\lfloor{n-1\over k-1}\right\rfloor - (k-1)\bigg(\left\lfloor{n-1\over k-1}\right\rfloor^2
+ \left\lfloor{n-1\over k-1}\right\rfloor\bigg)\adveq$$
for $2\leq k\leq n$,
as well as the bounds
\global\kpermsineq=\eqcount
$${(n-k+2)(n-1)\over k-1} - k + 1 \leq P_{nk} \leq {(n-k+2)(n-1)\over k-1} + k-3.\adveq$$
\noblueboxes

\proof The case $k=1$ counts the $n$ one-element sequences $(i)$ for $1\leq i\leq n$.
Let $2\leq k\leq n$; here the chief constraint is that both $a$ and $a+(k-1)r$ must belong to $[1,n]$.
We sum over the possible values of the step $r$. It must belong to one of the two discrete intervals
$\pm\big[1,\big((n-1)/(k-1)\big)\big]$,
and the number of valid choices for $a$ decreases as
$|r|$ increases. If $r$ is positive, then it is easy to see that $1\leq a \leq n-(k-1)r$ by the above constraint.
Similarly, if $r$ is negative, then $n + (k-1)r+1\leq a \leq n$.
In general, with $n$, $k$, and $r$ fixed, the number of possibilities
for $a$ is $n-(k-1)|r|$. Exploiting the symmetry of the two cases, we compute
$$\eqalign{
P_{nk} &= \sum_{r=1}^{\lfloor(n-1)/(k-1)\rfloor} 2\big(n-(k-1)r\big)\cr
&= 2n\left\lfloor{n-1\over k-1}\right\rfloor - 2(k-1)\sum_{r=1}^{\lfloor(n-1)/(k-1)\rfloor}r,\cr
}\adveq$$
yielding Equation~\refeq{\the\kpermseq}. Letting $\{x\}$ denote the fractional part of $x$ and using the identity
$\lfloor x\rfloor = x - \{x\}$, from Equation~\refeq{\the\kpermseq} we obtain
$$\eqalign{
P_{nk} &= 2n\bigg({n-1 \over k-1} - \left\{{n-1\over k-1}\right\}\bigg)\cr
&\qquad\qquad-(k-1)\bigg( \bigg( {n-1\over k-1}- \left\{{n-1\over k-1}\right\}\bigg)^2+ {n-1\over k-1}
-\left\{{n-1\over k-1}\right\}\bigg)\cr
&= \bigg({n-1 \over k-1}\bigg)\bigl(2n-(n-1)-(k-1)\bigr)\cr
&\qquad\qquad- \bigl(-2n+2(n-1)+(k-1)\bigr)\left\{{n-1\over k-1}\right\}
-(k-1) \left\{{n-1\over k-1}\right\}^2\cr
&={(n-k+2)(n-1)\over k-1} + (k-3)\bigg\{{n-1\over k-1}\bigg\} - (k-1)\bigg\{{n-1\over k-1}\bigg\}^2,
}\adveq$$
which, since both $\{x\}$ and $\{x\}^2$ are nonnegative and less than 1 for any $x$,
proves Equation~\refeq{\the\kpermsineq}.\slug

\newcount\kdperms
\kdperms=\thmcount
\proclaim Lemma \advthm. Let $P_{nk}$ be the number of $k$-orderings of $[1,n]$ that are an arithmetic
progression. The number of $k$-orderings of $[1,n]^d$ that form an arithmetic progression is
$$P_{nkd} = (P_{nk} + n)^d - n^d.\adveq$$

\proof Let $a,r\in [1,n]^d$. Writing $a = (a_1, a_2, \ldots, a_d)$ and $r= (r_1,r_2,\ldots,r_d)$, the arithmetic
progression
$$\big(a, a+r, a+2r, \ldots, a+(k-1)r\big)$$
has entries in $[1,n]^d$ if and only if the projection onto the $j$th coordinate
$$\big(a_j, a_j+r_j, a_j+2r_j, \ldots, a_j+(k-1)r_j\big)$$
is an arithmetic progression and has entries in $[1,n]^d$ for every $1\leq j\leq d$.
All but one of the projections may be a trivial progression, with $r=0$.
We did not count the trivial progressions in Lemma~\the\kperms\ because we were counting sequences without
repetition,
so by including the $n$ trivial progressions, the number of valid pairs $(a_j,r_j)$ increases to $P_{nk}+n$
for any $j$. Forming the product over all $j$, we arrive at $(P_{nk}+n)^d$, but we also counted $n^d$
undesirable trivial progressions in $\ZZ^d$. Subtracting them proves the lemma.\slug

We were quite careful with the calculations in these lemmas, but to prove our first main result, we will only
need the fact that $P_{nkd}\sim_n n^{2d}/(k-1)^d$. The bulk of the work is done in the next lemma, which
has been formulated more generally so that it may be reused in later sections. We will require the
gamma function $\Gamma(z) = \int_0^\infty x^{z-1}e^{-x}\,dx$, which, in this form,
is defined for all complex numbers with $\Re z > 0$ and for all nonnegative integers, $\Gamma(n+1) = n!$
(a classic treatment can be found in~\ref{15}).
Our main result relies on the fact that
%the gamma function is strictly increasing on
%the real interval $(\mu, \infty)$, where $\mu \approx 1.4616$, as well as the fact that the ratio
%$\Gamma(x+1)/\Gamma(x)$ approaches infinity
%as $x\to\infty$ along the real axis. Note that
for any $f:\NN\to \RR$, there exists a function $g:\RR\to\RR$,
analytic on all of $\RR$, such that $g(n) = f(n)$ for all $n\in \NN$. (This is a special case of a result of
basic complex analysis; see for example,~\ref{1}, page\ 197.)

\newcount\secondmoment
\secondmoment=\thmcount
\proclaim Lemma \advthm. Let $(A_n,Z_n)$ be a family of nonempty additive sets
indexed by the positive integers and let $P_n(k)$ denote the number of $k$-orderings of $A_n$
that are an arithmetic progression. Suppose that
\medskip
\item{i)} for any fixed $k\in \NN$, ${P_n(k+1)/ P_n(k)}\to 1$ as $n\to \infty$; and
\smallskip
\item{ii)} there exists a function $f:\NN\to\RR$ such that for all $n$, $P_n(k) < k!$ for all $k>f(n)$
and for any $O\big(f(n)\big)$ function $g(n)$, $P_n\big(g(n)\big)\to\infty$ as $n\to\infty$.
\medskip
\noindent Let $P^*_n(x)$ be the analytic continuation of $P_n$ to all of $\RR$ and
let $L_n$ denote the length of the longest arithmetic subsequence in an ordering of $A_n$
chosen uniformly at random.
There exists a function $\psi: \NN\to \RR$ that is $O\big(f(n)\big)$ and $\omega(1)$
such that for all positive integers $n$,
$${P^*_n(\psi(n))\over \Gamma(\psi(n)+1)} = 1\adveq$$
and $\pr\{\psi(n)-6 \leq L_n<\psi(n)+1\}\to 1$ as $n\to\infty$.

\proof We first prove the existence of the function $\psi$. Let $n\in\NN$ and note first that
$P^*_n(1) = |A_n|$, since this is the number of distinct base points.
If $|A_n|=1$, we set $\psi(n) = 1$. Otherwise,
we may consider the function $h_n(x) = P^*_n(x) - \Gamma(x+1)$,
which is analytic on the real interval $[1,\infty)$ and positive at $x=1$, since in this case $h_n(1) = |A_n| > 1$.
By hypothesis (ii), whenever $k>f(n)$ is an integer we must have $P^*_n(k) < \Gamma(k+1)$.
$h_n\big(f(n)+1\big) < 0$, so by the intermediate value
theorem, there exists a point $1<x^*<f(n)+1$ for which $h_n(x^*) = 0$ and we can set $\psi(n) = x^*$ (if there
is more than one choice for $x^*$, we pick the smallest one).
It is clear from this construction
that $\psi(n)$ must be $O\big(f(n)\big)$, and $\psi(n)$ must also be $\omega(1)$; otherwise
setting $g=\psi$ violates hypothesis (ii) by continuity of the function $\Gamma(x+1)$.

In the random ordering of $A_n$, there are ${|A_n| \choose k}$ subsequences of length $k$.
Each of the $k!{|A_n| \choose k} = |A_n|!/(|A_n|-k)!$ orderings of length $k$ has an equal chance of appearing
as a subsequence of the random ordering, so the
probability of a given subsequence of length $k$ being an arithmetic progression is $P_n(k) (|A_n|-k)! / |A_n|!$.
Let $\C(k)$ be the set of all subsequences of length $k$ in the original random sequence; we have
$|\C(k)| = {|A_n|\choose k}$. For a subsequence
$S\in\C(k)$, let $B_S$ be the event that $S$ is a $k$-term arithmetic progression.
Letting $N(n,k)$ denote the number of
$k$-term arithmetic progressions in the ordering of $A_n$, we have, by linearity of expectation,
$$\eqalign{
\ex\big\{N(n,k)\big\} &= \ex\Big\{\sum_{S\in \C(k)} \one_{B_S}\Big\} = \sum_{S\in \C(k)} \pr\{B_S\}\cr
&= {|A_n|\choose k} {P_n(k) (|A_n|-k)!\over |A_n|!} = {P_n(k)\over k!}.
}\adveq$$
By hypothesis (i), we see that
$${\ex\big\{N(n,k+1)\big\}\over \ex\big\{N(n,k)\big\}}
= {P_n(k+1)\over (k+1)!}\cdot{k!\over P_n(k)} \sim_n {1\over k+1}\adveq$$
and hence the expected value of $N(n,k)$ decreases with $k$. We constructed $\psi(n)$
so that $\ex\big\{N(n,\psi(n))\} = 1$, so a union bound now gives
$$\eqalign{
\pr\{L_n\geq \psi(n)+1\} & \leq \ex\big\{N(n,\psi(n)+1)\big\}\cr
&={\ex\big\{N(n,\psi(n)+1)\big\}\over \ex\big\{N(n,\psi(n))\big\}}\cr
&\sim_n {1\over \psi(n)+1},
}\adveq$$
which goes to zero as $n\to\infty$ since $\psi(n)\to \infty$.

On the other hand, we have
$${\ex\big\{N(n,\psi(n))\big\}\over \ex\big\{N(n,\psi(n)-1\big\}}\sim_n {1\over \psi(n)},\adveq$$
which implies that $\ex\big\{N(n,\psi(n)-1)\big\} \sim_n \psi(n)$. Continuing in this manner,
we see that $\ex\big\{N(n,\psi(n)-s)\big\} \sim_n \psi(n)^s$ and
$\lim_{n\to\infty} \ex\big\{N(n,\psi(n)-s)\big\}=\infty$ for all integers $1\leq s< \psi(n)$.
To turn this into a lower bound on $L_n$, one
must also analyse the behaviour of the second moment
$$\ex\big\{{N(n,k)}^2\big\}
= \ex\bigg\{ \Big(\sum_{S\in \C(k)}\one_{B_S}\Big) \Big(\sum_{T\in \C(k)} \one_{B_T}\Big)\bigg\}.
\adveq$$
For any $S$ in $\C(k)$ and $0\leq m\leq k$,
let $\C_{Sm}(k)\subseteq \C(k)$ denote the set of all subsequences $T$ of length $k$ that have exactly $m$
elements in common with $S$ and that are themselves progressions (i.e., such that $B_T$ holds).
Then, expanding the second moment to
$$\ex\big\{{N(n,k)}^2\big\} = \ex\Big\{\sum_{S\in \C(k)} \one_{B_S} \sum_{m=0}^k
\bigl|\C_{S_m}(k)\bigr|\Big\},\adveq$$
we approach the inner summation by considering each set $\C_{Sm}(k)$ separately.
Assuming both $S$ and $T$ are $k$-term arithmetic progressions,
for $1\leq m\leq k$, we note that if $T$ and $S$ both contain $m$ elements, then those $m$ elements must form
a subprogression of both $T$ and $S$. There are no more than $k^2$ ways to form a subprogression of length
$m$ in $S$, since there are $k$ choices for the first element of the subprogression and $\leq k$ choices for
the second element; these two choices determine the subprogression entirely.
Now, given a specific subprogression $R$ of length $m$, we count the number of possibilities for $T$. The
first element of $R$ is the $j$th element of $T$ for some $1\leq j\leq k$,
so there are no more than $k$ possibilities for its position in $T$.
There are also $\leq k$ choices for the second element of $R$, and after this, $T$ is completely defined.
Thus $|\C_{Sm}(k)| \leq k^4$ for all $1\leq m\leq k$.
(This is a substantial overcount for $m>1$; when $m=k$ there is in fact only one possibility for $T$.)
Next, we consider the set $\C_{S0}(k)$ of all sequences that have no elements in common with $S$,
which has size ${|A_n| - k\choose k}$.
Each of these subsequences is
a progression with probability no more than
$${P_n(k)\over{|A_n|-k\choose k}k!} = {P_n(k)(|A_n| - 2k)!\over(|A_n| - k)!},\adveq$$
so the expected number of progressions
among these subsequences is at most $P_n(k)/k!$. Putting these facts together, we find that
$$\ex\big\{{N(n,k)}^2\big\} \leq {P_n(k)\over k!}\bigg( {P_n(k)\over k!} + \sum_{m=1}^k k^4 \bigg)
= {P_n(k)\over k!}\bigg( {P_n(k)\over k!} + k^5\bigg),\adveq$$
\noblueboxes meaning that $\var\big\{N(n,k)\big\} \leq  P_n(k)k^5/k!$.
We have, by the Chung-Erd\H os inequality~\ref{4},
$$\pr\{L_n< k\} = \pr\big\{{N(n,k)}=0\big\} \leq {\var\big\{N(n,k)\big\}\over \ex\big\{N(n,k)\big\}^2}
= {k^5\over P_n(k)/k!}.$$
Plugging in $k= \psi(n)-6$, we have
$$\pr\{L_n < \psi(n)-6\} \leq {\var\big\{N(n,\psi(n)-6)\big\}\over \ex\big\{N(n,\psi(n)-6)\big\}^2}
\sim_n {\big(\psi(n)-6\big)^5\over \psi(n)^6} \to 0.$$
We have shown that $\pr\big\{\psi(n)-6\leq L_n < \psi(n) +1\big\}\to 1$, which is what we wanted.\slug

Note that this lemma does not use the commutativity of the groups $Z_n$; in the final section of this paper
we briefly discuss possible noncommutative generalisations.
Our first main theorem is a direct consequence of the previous lemma.

\newcount\torsionfreethm
\torsionfreethm=\thmcount
\proclaim Theorem \advthm. For positive integers $n$ and $d$,
let $L_n$ denote the longest arithmetic subsequence in an ordering of $[1,n]^d$,
chosen uniformly at random. There exists a function $\psi(n,d)$ with $\psi(n,d)\sim_n 2d\log n/\log\log n$ such
that the probability that $\psi(n,d)-6\leq L_n <\psi(n,d)+1$
tends to $1$ as $n$ approaches infinity.
% Required on every page that has blue boxes
\special{ps:[/pdfm { /big_fat_array exch def big_fat_array 1 get 0
0 put big_fat_array 1 get 1 0 put big_fat_array 1 get 2 0 put big_fat_array pdfmnew } def}

\proof Set $Z_n = \ZZ^d$ for all $n$, $A_n = [1,n]^d$ and let $P^*_n(k)$ be the analytic continuation of
$P_{nkd}$ as defined in Lemma~\the\kdperms. A simple computation shows that for fixed $k$ and $d$,
$P_{(n+1)kd}/P_{nkd}\to 1$ as $n\to\infty$, and setting $f(n) = n$, we have $P_{nkd} = 0 < k!$ for all $k>f(n)$.
By our earlier observation that $P_{nkd}\sim_n n^{2d}/(k-1)^d$, substituting
any $O(n)$ function for $k$ we have $\lim_{n\to\infty} P_{nkd}= \infty$.
Thus we may apply Lemma~{\the\secondmoment} to obtain a function $\psi(n,d)$ satisfying
$${P_{n(\psi(n,d))d}\over \psi(n,d)!} = 1,$$
for all $n\in \NN$. The fact that $P_{nkd}\sim_n n^{2d}/(k-1)^d$ implies that
$\psi(n,d)!\sim_n n^{2d}/\big(\psi(n,d)-1\big)^d$, and Stirling's formula gives
$\psi(n,d)\sim_n 2d\log n/\log\log n$. The rest of the theorem follows from the previous lemma, and
{\it a fortiori} we have $L_n\sim_n 2d\log n/\log\log n$.\slug

Because the set of arithmetic progressions is invariant under translation and dilation, as a corollary
in the integer case we can replace $[1,n]$ with the entries of any arithmetic progression of length $n$.

\proclaim Corollary \advthm. Let $S$ be an arithmetic progression of length $n$ in $\ZZ$, i.e., $S$ is of the form
$$S = \big(a,a+r,a+2r,\ldots, a+(n-1)r\big)$$
for some $a,r\in \ZZ$, $r\neq 0$. Randomly permute the entries of $S$ and let $L_n$ be the length of the longest
arithmetic progression in the shuffled sequence. If $\psi$ is the function given by Theorem~\the\torsionfreethm,
then $\psi(n)-6\leq L_n< \psi(n)+1$
with probability tending to 1 as the length $n$ of the arithmetic progression tends to infinity.
\slug

\medskip
\hldest{xyz}{}{countexact}
\boldlabel Counting exact values.
For $d=1$ and small values of $n$, one can easily
generate the exact distribution of $L_n$ using a computer. Let
$f_n(k)$ denote the number of permutations of $[1,n]$ whose longest embedded arithmetic progression has
length $k$; so $\pr\{L_n = k\} = f_n(k)/n!$. Table 1 contains the values of $f_n(k)$ for $n\leq 10$. No general
formula for $f_n(k)$ has been found, but the case $k=2$ (these are known as {\it $3$-free permutations})
has received a fair amount of attention in the literature.
A. Sharma found the explicit upper bound $f_n(2)\leq (2.7)^n/21$, valid for $n\geq 11$, as well as the
asymptotic result $f_n(2) = \omega_n(2^nn^k)$, which holds for any choice of $k$~\ref{12}.
The lower bound was improved
to $(1/2)c^n$ for $n\geq 8$, where $c=2132^{1/10}\approx 2.152$,
by T. D. LeSaulnier and S. Vijay~\ref{11}.
More recently, B. Correll, Jr.\ and
R. W. Ho presented an efficient dynamic programming algorithm to count $3$-free permutations and used these
results to inductively refine previously-known bounds~\ref{3}.

\topinsert
$$\vbox{
\centerline{\smallheader Table 1}
\medskip
\centerline{\eightpoint PERMUTATIONS OF $[1,n]$ WHOSE LONGEST PROGRESSION HAS LENGTH $k$}
}
$$
$$\centerline{\vbox{
\eightpoint
\hrule
\medskip
\tabskip=.5em plus.2em minus .5em
\halign{
   $\hfil#$   &  $\hfil#$ & $\hfil#$ & $\hfil#$ & $\hfil#$ & $\hfil#$ & $\hfil#$ & $\hfil#$ &
   $\hfil#$ & $\hfil#$ & $\hfil#$ & $\hfil#$ \cr
   n   &  f_n(1) & f_n(2) & f_n(3) & f_n(4) & f_n(5) & f_n(6) & f_n(7) & f_n(8) & f_n(9) & f_n(10) \cr
   \noalign{\medskip}
   \noalign{\hrule}
   \noalign{\medskip}
1 & 1 \cr
2 & 0 & 2 \cr
3 & 0 & 4 & 2 \cr
4 & 0 & 10 & 12 & 2 \cr
5 & 0 & 20 & 82 & 16 & 2 \cr
6 & 0 & 48 & 516 & 134 & 20 & 2 \cr
7 & 0 & 104 & 3232 & 1480 & 198 & 24 & 2 \cr
8 & 0 & 282 & 21984 & 15702 & 2048 & 274 & 28 & 2 \cr
9 & 0 & 496 & 168368 & 162368 & 28048 & 3204 & 362 & 32 & 2 \cr
10 & 0 & 1066 & 1306404 & 1902496 & 374194 & 39420 & 4720 & 462 & 36 & 2 \cr
 \noalign{\medskip}
 \noalign{\hrule}
    }
}}$$
\endinsert

\advsect Finite additive groups
\hldest{xyz}{}{sec3}

The aim of this section is to extend our result for torsion-free groups to finite abelian groups.
Thus, for a finite abelian group $Z$, we are now studying the case $A = Z$. The extension
proceeds smoothly in the case of cyclic groups, which in turn provides some information for arbitrary
finite additive groups via the structure theorem.

\medskip
\hldest{xyz}{}{cyclic}
\boldlabel Cyclic groups. We shall prove that the case in which $Z= \zed n$ and $A = Z$ has the same
asymptotics as the case $A = [1,n]\subseteq \ZZ$.
Before we proceed, we recall that
Euler's {\it totient function} $\tot(n)$ counts the numbers coprime to and less than a given $n$, that is,
$$\tot(n) = \big|\{m\in [1,n] : \gcd(m,n) = 1\}\big|.$$
A basic result of group theory states that for $1\leq d< n$, the number of elements of order $d$ in $\zed n$ is
$\tot(d)$. We will write $d\divides n$ to indicate that $d$ divides $n$.
As before, we count progressions in $k$-orderings of $A$, starting first with the case that $Z = \zed n$,
the cyclic group of order $n$.

\newcount\kpermscyclic
\kpermscyclic=\thmcount
\newcount\kpermscycliceq
\newcount\kpermscyclicineq
\proclaim Lemma \advthm. Let $Q_{nk}$ denote the number of $k$-orderings of $\ZZ/n\ZZ$
that form an arithmetic progression. For every $n$, $Q_{n1} = n$ and $Q_{nk} = 0$ for all $k>n$; for
$2\leq k\leq n$, we have
\global\kpermscycliceq=\eqcount
$$Q_{nk} = n\Big(n-\sum_{j=1}^{k-1} \indic{j\divides n}\tot(j)\Big).\adveq$$
In particular, $Q_{nn} = n\tot(n)$ and we also have the bounds
\global\kpermscyclicineq=\eqcount
$$n^2 - {3\over \pi^2}n(k-1)^2 - o_k(nk^2)\leq Q_{nk} \leq n(n-1).\adveq$$
If $n$ is prime, then $Q_{nk} = n(n-1)$ for all $k$.

\proof The cases $k=1$ and $k>n$ are obvious, so let $2\leq k\leq n$.
By the $n$-fold symmetry of $\zed n$, it suffices to count the number of $k$-term arithmetic progressions
with base point
$0$, and then multiply the result by $n$. For any $r\in \zed n\setminus\{0\}$, the sequence
$$\big(0, r, 2r, \ldots, (k-1)r\big)$$
is a valid $k$-ordering if and only if the element $0$ does not appear twice.
This is equivalent to the condition that the order of $r$ be at least $k$. Since the order of a group
element must divide $n$, to calculate the number of possible steps $r$ we start with $n$ and subtract
the number of elements of order $j$ for $1\leq j\leq k$ where $j\divides n$; this equals $\tot(j)$ for all $j<n$.
Multiplying this by the $n$ possible base points yields the formula Equation~\refeq{\the\kpermscycliceq}.

The fact that $Q_{nn} = n\tot(n)$ follows from the observation that the summation equals $n-\tot(n)$ when $k=n$.
The upper bound in
Equation~\refeq{\the\kpermscyclicineq}\ is clear because the summation is at least $1$ (for $j=1$); if
$n$ is prime, then this bound is met with equality. To prove the lower bound, we use the
asymptotic formula
$$\sum_{j=1}^{k} \tot(k) = {3\over \pi^2}k^2 +
  O_k\big(k(\log k)^{2/3}(\log\log k)^{2/3}\big),\adveq$$
\noblueboxes
whose proof can be found in \ref{14}. It yields
$$\eqalign{
Q_{nk} &=n\Big(n-\sum_{j=1}^{k-1} \indic{j\divides n}\tot(j)\Big)\cr
&\geq n\Big(n-\sum_{j=1}^{k-1}\tot(j)\Big)\cr
&= n^2 - {3\over \pi^2}n(k-1)^2 - o_k(nk^2).\noskipslug
}\adveq$$

Now we can compute the length of the longest arithmetic subsequence in an ordering of a cyclic
group, which is similar to the case $d=1$ in Theorem~\the\torsionfreethm.

\newcount\cyclicthm
\cyclicthm=\thmcount
\proclaim Theorem \advthm. For a positive integer $n$, let $L_n$ denote the length of the longest arithmetic
subsequence of an ordering of $\ZZ/n\ZZ$, chosen uniformly at random.
There exists a function $\chi(n)$ with
$\chi(n)\sim_n 2\log n/\log\log n$ such that $\chi(n)-6\leq L_n<\chi(n)+1$
with probability tending to 1 as $n$ tends to infinity.

\proof
Let $Q_{nk}$ be as in
Lemma~\the\kpermscyclic\ and
let $Q^*_n(x)$ be an analytic function on all of $\RR$ such that $Q^*_n(k) = Q_{nk}$ for all positive integers $k$.
For a fixed $k$, we have
$${(n+1)^2 - O_k\big((n+1)k^2\big)\over n(n-1)}\leq{Q_{(n+1)k}\over Q_{nk}}\leq{(n+1)n\over n^2-O_k(nk^2)}.\adveq$$
Taking $n$ to infinity in both bounds tells us that $Q_{(n+1)k}/Q_{nk}\to 1$.
Let $f(n) = n^{1/3}+C$, where $C$ is chosen such that $f(n)! > n(n-1)$ for all $n\in \NN$.
So $Q_{nk} \leq n(n-1) < f(n)! < k!$ for any $k>f(n)$ and letting $g(n)$ be any $O(n^{1/3})$ function,
we have, by the previous lemma, $Q_{ng(n)} \geq n^2 - O(ng(n)^2) = n^2 - o(n^2)\to \infty$, so we can
apply Lemma~{\the\secondmoment} with $A_n = Z_n = \zed n$ for all $n$ to obtain a function $\chi$ satisfying
$${Q^*_n(\chi(n))\over\Gamma(\chi(n)+1)} = 1.$$
Since $\chi(n)! = Q_{n\chi(n)} = n^2 - o(n^2)$ for all $n\in \NN$, another application of Stirling's formula
gives $\chi(n)\sim_n 2\log n/\log\log n$,
and we have $\chi(n)-6\leq L_n< \chi(n)+1$ with
probability tending to $1$.\slug

Note that although the function $\psi(n,1)$ obtained by applying Theorem~\the\torsionfreethm\ to $d=1$
has the same asymptotics as $\chi(n)$ from Theorem~\the\cyclicthm, we have $\psi(n,1) < \chi(n)$. This
follows from the functional equations defining $\psi$ and $\chi$ and also makes sense intuitively, since
it is easier, in some sense, to form a progression if one is allowed to loop around the edge of the interval,
which is possible in the cyclic
case. As a concrete example, $(0,2,6,1,3,5,4)$ has the $4$-term arithmetic subsequence $(0,6,5,4)$ when the
sequence is regarded as an ordering of the group $\ZZ/7\ZZ$ (the base point is $0$ and the step size is $6$),
but it has no arithmetic subsequence of length $4$ when regarded simply as an ordering of the discrete
interval $[0,6]$.

\medbreak
\hldest{xyz}{}{exactval}
\boldlabel Exact values.
Table 2 contains exact values of $g_n(k) = n!\pr\{L_n = k\}$ for small $n$, and comparing
Tables 1 and 2, it is clear that $L_n$ is expected to be longer when $A_n = \zed n$ than when $A_n = [1,n]$.
A closer look at Table 2
reveals certain curiosities not present in the torsion-free case. Firstly, the sequence $\ex\{L_n\}$
is not strictly increasing; for instance, we have $\ex\{L_7\} = 4.25$ and $\ex\{L_8\} \approx 4.136$. The same
can be said for $L_{11}$ and $L_{12}$, and this phenomenon can be attributed to the fact that $Q_{nk}$ is
greater, in proportion to $n$, when $n$ is prime. Next, we see that $\pr\{L_n = n\} = n\tot(n)$, which follows
directly from Lemma~\the\kpermscyclic. Lastly, we note that an ordering
of $\ZZ/n\ZZ$ can be 3-free only if $n$ is a power of $2$. In fact, we can derive a simple explicit formula
for the number of $3$-free orderings of $\zed n$.

\topinsert
$$\vbox{
\centerline{\smallheader Table 2}
\medskip
\centerline{\eightpoint PERMUTATIONS OF $\zed n$ WHOSE LONGEST PROGRESSION HAS LENGTH $k$}
}
$$
$$\centerline{\vbox{
\eightpoint
\hrule
\medskip
\tabskip=.7em plus.2em minus .5em
\halign{
   $\hfil#$   &  $\hfil#$ & $\hfil#$ & $\hfil#$ & $\hfil#$ & $\hfil#$ & $\hfil#$ & $\hfil#$ &
   $\hfil#$ & $\hfil#$ & $\hfil#$ \cr
   n   &  g_n(1) & g_n(2) & g_n(3) & g_n(4) & g_n(5) & g_n(6) & g_n(7) & g_n(8) & g_n(9) & g_n(10) \cr
   \noalign{\medskip}
   \noalign{\hrule}
   \noalign{\medskip}
1 & 1 \cr
2 & 0 & 2 \cr
3 & 0 & 0 & 6 \cr
4 & 0 & 8 & 8 & 8 \cr
5 & 0 & 0 & 40 & 60 & 20 \cr
6 & 0 & 0 & 468 & 192 & 48 & 12 \cr
7 & 0 & 0 & 462 & 3150 & 1176 & 210 & 42 \cr
8 & 0 & 128 & 4192 & 27872 & 6592 & 1312 & 192 & 32 \cr
9 & 0 & 0 & 57402 & 182790 & 99630 & 19656 & 2970 & 378 & 54 \cr
10 & 0 & 0 & 67440 & 1795320 & 1594640 & 146200 & 22000 & 2840 & 320 & 40 \cr
 \noalign{\medskip}
 \noalign{\hrule}
    }
}}$$
\endinsert

\proclaim Theorem \advthm. Let $n\geq 1$ be an integer.
The number $g_n(2)$ of orderings of $\zed n$ that do not contain any arithmetic
subsequence of length 3 equals $2^{n-1}$ if $n=2^m$ for some $m\geq 1$, and is zero otherwise. An ordering
of $\zed{2^m}$ that contains no progression of length $3$ consists of $2^{m-1}$ elements of the same parity,
followed by the $2^{m-1}$ elements of the opposite parity.

\proof In this proof all arithmetic operations are taken modulo $n$. Obviously $g_1(2) = 0$, so we begin
by supposing that $n = p$ is an odd prime and fixing an arbitrary ordering
of $\zed p$. Let $z_1$ and $z_2$ denote the first two elements of the sequence and consider $r=z_2-z_1$.
Since $p$ is odd, $z_2+r\neq z_1$ and must therefore come later in the sequence. So there is an embedded
progression of length $3$. More generally, suppose that $n\geq 3$ is not a power of $2$ and let $p$ be an odd
prime that divides $n$. Note that the elements $0, n/p, 2n/p, \ldots, (p-1)n/p$ must appear in some order
in the sequence, and this sequence contains a $3$-term progression if and only if the same ordering, with
each element divided by $n/p$, is a $3$-term progression. But this is an ordering of $0,1,2,\ldots,p-1$,
and we already showed that every ordering of $\zed p$ contains a $3$-term arithmetic progression when $p$
is an odd prime.

Now we handle the case in which $n=2^m$, by induction on $m$.
For the case $m=1$, both orderings of $\ZZ/2\ZZ$ are $3$-free, and both can be split into an odd half
and an even half. Now suppose that there are $2^{n/2-1}$
orderings of $\zed{(n/2)}$ that are $3$-free. For any pair $(S,T)$ of such orderings, with
$S = (s_1, s_2,\ldots,s_{n/2})$ and $T=(t_1,t_2,\ldots,t_{n/2})$, note that
$$(2s_1,2s_2,\ldots,2s_{n/2}, 2t_1+1, 2t_2+1, \ldots, 2t_1+1)$$
and
$$(2t_1+1,2t_2+1,\ldots,2t_{n/2}+1, 2s_1,2s_2,\ldots,2s_{n/2})$$
are two orderings of $\zed n$ that are $3$-free and uniquely determined by the pair $(S,T)$. This means there
are at least $2\cdot (2^{n/2-1})^2 = 2^{n-1}$ orderings of $\zed n$ with no $3$-term arithmetic progressions.
To see that we have, in fact, counted all of them, let an ordering of $\zed n$ be given that is $3$-free;
we aim to show that it is of one of the two forms prescribed above. Let $S$ be
the subsequence of all even elements, and let $T$ be the subsequence of all odd elements. Note that neither
of these two subsequences can contain a progression of length $3$, so these sequences $S$ and $T$
were included in the count above.
The last thing to show is that all the even elements are on one side of the ordering and all of the odd
elements are on the other side. To do this, fix an odd element $z$ and consider $z+1$ and $z-1$. These elements
are distinct (since $n$ is a power of $2$) and both even,
but if one is to the left of $z$ and one is to the right, then a $3$-term progression appears.
So both $z+1$ and $z-1$ are on the same side of $z$ (without loss of generality, assume it is the left side).
Next, consider $z+3$ and $z-3$. For the same reason
as before, they must both be on the same side of $z$, but if they are to the right of $z$,
then we would have one of the arithmetic subsequences $(z-1, z+1, z+3)$ or $(z+1, z-1, z-3)$, depending
on the order of $z+1$ and $z-1$. An analogous argument applies if $z+1$ and $z-1$ had been to the right of $z$.
So $z+3$ and $z-3$ are on the same side of $z$ as $z+1$ and $z-1$. Continuing
in this manner for $z\pm 5, z\pm 7,\ldots$, we find that all of the even elements are on the same side of
all of the odd elements, meaning
that the sequences of the form given in the theorem statement are the only $3$-free orderings. The induction is
complete and the theorem is proved.\slug

\medskip
\hldest{xyz}{}{arbfinite}
\boldlabel Arbitrary finite abelian groups.
Because the groups of order $n$ may have wildly different structures, we cannot hope to give a sensible
asymptotic result that is valid for all finite additive groups. For example, both $\zed{2^d}$ and $(\zed 2)^d$
have order $2^d$, but Theorem~\the\cyclicthm\ tells us that a random ordering of the first group will
have an arithmetic subsequence of length
asymptotic to $2\log(2^d)/\log\log(2^d)$ with probability tending to $1$,
while an ordering of the second group cannot have a 3-term arithmetic subsequence.
Every finite abelian group $Z$ admits a unique {\it invariant factor decomposition}
$$Z = \zed{n_1}\times \zed{n_2}\times\cdots\times\zed{n_d}$$
where $n_1\divides n_2$, $n_2\divides n_3$, and so on (the factors increase in size from left to right).
In other words, $\zed{n_d}$ is the largest cyclic
subgroup of $Z$ and $Z$ cannot be written as a product of fewer than $d$ cyclic groups. Of course, we have
$n = n_1n_2\cdots n_d$. With this decomposition,
we can generalise the bounds obtained in Lemma~\the\kpermscyclic\ to all finite additive groups.

\proclaim Lemma \advthm. Let $Z$ be a group of order $n$ whose invariant factor decomposition is
$$Z = \zed{n_1}\times \zed{n_2}\times\cdots\times\zed{n_d}.$$
Let $Q_{nk}$ be as in Lemma~\the\kpermscyclic\ and for $2\leq k\leq n_d$, let $Q_k(Z)$ be the number of
$k$-orderings of $Z$ that are an arithmetic progression; so $Q_{nk} = Q_k(\zed n)$. Let $j$ be the smallest
index for which $k\leq n_j$. Then
$$n\prod_{i=j}^d n_i - O_k(nk^{2d}) \leq Q_{k}(Z) \leq n\prod_{i=j}^d n_i - n.\adveq$$

\proof For the upper bound, we cannot simply take a product over the $Q_{n_ik}$,
because this does not count progressions
whose projection onto some factor is not a valid $k$-ordering. For example, $\big((0,4), (2,6), (0,0)\big)$ is
a valid progression in $\zed 4\times\zed 8$, but its projection onto $\zed 4$ is $(0,2,0)$, which is not a valid
3-ordering. Instead, we form the product of trivial progressions in the first $j-1$ factors with
the number ${n_i}^2$ of starting pairs in the last $d-j+1$ factors, and then subtract the $n$ progressions
that are trivial in $Z$, obtaining
$$Q_{nk} \leq \Big(\prod_{i=1}^{j-1} n_i\Big)\Big(\prod_{i=j}^d {n_i}^2\Big) - n = n\prod_{i=j}^d n_i - n.\adveq$$
As for the lower bound,
we can undercount the number of progressions in $Z$ by counting the product of trivial progressions
in the first $j-1$ factors and valid nontrivial progressions in the remaining $d-j+1$ factors. This gives~us
$$\eqalign{
Q_k(Z) &\geq \Big(\prod_{i=1}^{j-1} n_i\Big)\Big(\prod_{i=j}^d Q_{n_ik}\Big)\cr
&\geq \Big(\prod_{i=1}^{j-1} n_i\Big)\Big(\prod_{i=j}^d \big({n_i}^2 - 3{n_i}(k-1)^2/\pi^2 -
o_k(n_ik^2)\big)\Big)\cr
&\geq n\prod_{i=j}^d n_i - O_k(nk^{2d}).\noskipslug\cr
}\adveq$$

Thus for a finite abelian group
$Z$ and a given $k$, one can compute the decomposition of $Z$ into $d$ invariant factors $\zed{n_1},\ldots,
\zed{n_d}$ and
determine the smallest index $j$ for which $k\leq n_j$. Then if $N_k$ denotes the number of $k$-term arithmetic
progressions that appear as subsequences in a random ordering of $Z$, we have
$${nn_jn_{j+1}\cdots n_d -n\over k!} - O_k(nk^{2d}) \leq \ex\{N_k\} \leq {nn_jn_{j+1}\cdots n_d -n\over k!}.\adveq$$
\goodbreak
This fact can be used to construct sequences of abelian groups $Z_n$ that can be fed into
Lemma~\the\secondmoment\ to obtain asymptotic information about the length of the maximal arithmetic subsequence.

\medskip
\hldest{xyz}{}{elempgroups}
\boldlabel Elementary {\mathbold p}-groups. We end this section by outlining an
an explicit example in which $n_1 = n_2 =\cdots = n_d = p$, a prime number.
An {\it elementary $p$-group} is an abelian group isomorphic
to $(\zed p)^d$ for some prime $p$ and integer $d$.
Thus we have a sequence of additive groups $Z_p = (\zed p)^d$ indexed by the primes.
Since $Q_{pk} = p(p-1)$ for every $2\leq k\leq p$,
by a calculation similar to that used in the proof of Lemma~\the\kdperms, the number of $k$-orderings of
$(\zed p)^d$ that are an arithmetic progression is $(Q_{pk} + p)^d - p^d = p^{2d} - p^d$. Letting $\tau(p,d)$
be the function from $\RR$ to $\RR$ satisfying
$${p^{2d} - p^d\over \Gamma(\tau(p,d)+1)} = 1,\adveq$$
we have $\tau(p,d) \sim_p 2d\log p/\log\log p$ and Lemma~\the\secondmoment\ tells us that
$\tau(p,d)-6\leq L_n\tau(p,d)+1$
with probability tending to one as $p$ marches towards infinity along the primes.
(Strictly speaking, one would actually have to modify the statement and hypotheses
of Lemma~\the\secondmoment\ to allow for
sequences of groups $A_n$ indexed by any monotone increasing sequence of positive integers.)
The same asymptotic result holds if the primes $p$ are replaced by arbitrary integers $n$,
but in this case the function $\tau(n,d)$ is not monotone with respect to $n$. So although we still have
$\tau(n,d) \sim_n 2d\log n/\log\log n$, the values of $\tau(n,d)$
will fluctuate as $n$ approaches infinity, attaining local maxima at prime values of~$n$.

\advsect Further work
\hldest{xyz}{}{sec4}

We will now discuss some related problems that have yet to be resolved, as well as possible noncommutative
generalisations.
As mentioned in the introduction, it is still unknown whether permutations of the positive integers
must necessarily contain arithmetic progressions of length 4. Furthermore, the functions $f_n(k)$ and $g_n(k)$
remain largely unstudied, except in the case $k=2$. One might study the number of $(k+1)$-free orderings
of $[1,n]$ or $\ZZ/n\ZZ$, which corresponds to the partial sums
$$\sum_{j=1}^k f_n(j)\quad\hbox{and}\quad\sum_{j=1}^k g_n(j),$$
respectively.

We end off by briefly considering what happens when we remove the condition that the underlying group
be abelian. We will denote a non-abelian group by $G$ and use multiplicative notation;
thus for $r\in G$ we have, for instance, $r^2 = rr$ and $r^{-3} = r^{-1}r^{-1}r^{-1}$. A {\it left
progression} is a sequence $(a,ra,r^2a, \ldots, r^{k-1}a)$, where $a,r\in G$ and $k\in \ZZ$. With the same
definitions for $a$, $r$, and $k$, a {\it right progression} is a sequence $(a,ar, ar^2,\ldots,
ar^{k-1})$. Note that the length of the longest left progression in a sequence is not the same,
in general, as the length of the longest right progression. For example, in the free group $F_2$,
with generators $a$ and $b$, the sequence $(a,ba,b^2a)$ contains a left progression of length 3
but no right progression of length 3.

However, in the case that $A\subseteq G$ is closed under inverses, the number of $k$-orderings of $A$
that are left progressions is the same as the number of $k$-orderings of $A$ that are right
progressions, because if $S(A,k)$ is the set of all $k$-orderings of $A$, there is a bijection
$f:S(A,k)\to S(A,k)$ that maps $(s_1,s_2,\ldots,s_k)$ to $({s_1}^{-1}, {s_2}^{-1},\ldots,{s_k}^{-1})$,
and a sequence $T\in S(A,k)$ is a left progression if and only if $f(T)$ is a right
progression. In particular, if every $k$-ordering of $A$ is just as likely to arise as a subsequence of
some ordering of $A$, then the expected number of left progressions of length $k$ in a random ordering
is the same as the expected number of right progressions.

It would be interesting to compute asymptotic formulas for $L_n$, akin to the ones we found for $\ZZ^d$, $\zed n$,
and $(\zed p)^d$, for families of non-abelian groups, such as the free groups $F_d$, the symmetric groups
$\frak S_n$, or the dihedral groups $D_n$. In an infinite group such as $F_d$, one must select a
finite subset $A$. A somewhat natural choice for the free group is the set of all reduced words
with length at most $n$, but perhaps a more faithful analogue of $[1,n]^d\subseteq \ZZ^d$ is the set of all words
$w\in F_d$
such that $f(w)\in [-n,n]^d$ where $f$ is the canonical homomorphism from $F_d$ to $F_d/[F_d,F_d]\cong \ZZ^d$.
In any case, both of these subsets of $F_d$ are closed under inverses, so
it suffices to study, say, the longest right progression. The second moment method we used in
Lemma~{\the\secondmoment} did not rely
on the group being abelian, but we anticipate that more sophisticated counting methods will be required to
count the expected number of progressions
in subsequences of orderings when dealing with nonabelian groups.

\hldest{xyz}{}{acks}
\section Acknowledgements

The authors are indebted to Luc Devroye for his numerous corrections and insightful comments.
We would also like to thank our homies Jaya Bonelli, Anna Brandenberger, Joseph Dahdah,
Jad Hamdan, Denali Relles, and Jonah Saks
for their encouragement and technical feedback. Lastly, we thank the anonymous referee for spotting some errors
in the preliminary draft, and for suggesting many improvements to the readablity of the paper and
the clarity of the proofs.
Both authors are supported by the Natural Sciences and Engineering Research Council of Canada.

\goodbreak
\hldest{xyz}{}{refs}
\section References

\parskip=0pt
\hyphenpenalty=-1000 \pretolerance=-1 \tolerance=1000
\doublehyphendemerits=-100000 \finalhyphendemerits=-100000
\frenchspacing
\def\bref#1{[#1]}
\def\beginref{\noindent
}
\def\endref{\medskip}
\vskip\parskip

\beginref \parindent=20pt\item{\bref{1}}
\hldest{xyz}{}{bib1}%
Lars Valerian Ahlfors,
{\sl Complex Analysis}
(New York:
McGraw-Hill,
1953).
\endref

\beginref \parindent=20pt\item{\bref{2}}
\hldest{xyz}{}{bib2}%
Itai Benjamini,
Ariel Yadin,
and Ofer Zeitouni,
``Maximal arithmetic progressions in random subsets,''
{\sl Electronic Communications in Probability}\/
{\bf 12}
(2007),
365--376.
\endref

\beginref \parindent=20pt\item{\bref{3}}
\hldest{xyz}{}{bib3}%
Bill Correll, Jr.~$\!\!$
and Randy Wesley Ho,
``A note on 3-free permutations,''
{\sl Integers: Electronic Journal of Combinatorial Number Theory}\/
{\bf 17}
(2017),
A55, 10 pp.
\endref

\beginref \parindent=20pt\item{\bref{4}}
\hldest{xyz}{}{bib4}%
Kai Lai Chung
and Paul Erd\H{o}s,
``On the application of the Borel-Cantelli lemma,''
{\sl Transactions of the American Mathematical Society}\/
{\bf 72}
(1952),
179--186.
\endref

\beginref \parindent=20pt\item{\bref{5}}
\hldest{xyz}{}{bib5}%
James Avery Davis,
Roger Charles Entringer,
Ronald Lewis Graham,
and Gustavus James Simmons,
``On permutations containing no long arithmetic progressions,''
{\sl Acta Arithmetica}\/
{\bf 34}
(1977),
81--90.
\endref

\beginref \parindent=20pt\item{\bref{6}}
\hldest{xyz}{}{bib6}%
Noam David Elkies
and Ashvin Anand Swaminathan,
``Permutations that destroy arithmetic progressions in elementary $p$-groups,''
{\sl Electronic Journal of Combinatorics}\/
{\bf 24}
(2017),
P1.20.
\endref

\beginref \parindent=20pt\item{\bref{7}}
\hldest{xyz}{}{bib7}%
Roger Charles Entringer
and Douglas Elmer Jackson,
``Elementary Problem 2440,''
{\sl American Mathematical Monthly}\/
{\bf 80}
(1973),
1058.
\endref

\beginref \parindent=20pt\item{\bref{8}}
\hldest{xyz}{}{bib8}%
Alan Frieze,
``On the length of the longest monotone subsequence in a random permutation,''
{\sl Annals of Applied Probability}\/
{\bf 1}
(1991),
301--305.
\endref

\beginref \parindent=20pt\item{\bref{9}}
\hldest{xyz}{}{bib9}%
Jesse Geneson,
``Forbidden arithmetic progressions in permutations of subsets of the integers,''
{\sl Discrete Mathematics}\/
{\bf 342}
(2019),
1489--1491.
\endref

\beginref \parindent=20pt\item{\bref{10}}
\hldest{xyz}{}{bib10}%
Peter Hegarty,
``Permutations avoiding arithmetic patterns,''
{\sl Electronic Journal of Combinatorics}\/
{\bf 11}
(2004),
R39.
\endref

\beginref \parindent=20pt\item{\bref{11}}
\hldest{xyz}{}{bib11}%
Timothy Dale LeSaulnier
and Sujith Vijay,
``On permutations avoiding arithmetic progressions,''
{\sl Discrete Mathematics}\/
{\bf 311}
(2011),
205--207.
\endref

\beginref \parindent=20pt\item{\bref{12}}
\hldest{xyz}{}{bib12}%
Arun Sharma,
``Enumerating permutations that avoid three-term arithmetic progressions,''
{\sl Electronic Journal of Combinatorics}\/
{\bf 16}
(2009),
R63.
\endref

\beginref \parindent=20pt\item{\bref{13}}
\hldest{xyz}{}{bib13}%
Endre Szemer\'edi,
``On sets of integers containing no $k$ elements in arithmetic progression,''
{\sl Acta Arithmetica}\/
{\bf 27}
(1975),
199--245.
\endref

\beginref \parindent=20pt\item{\bref{14}}
\hldest{xyz}{}{bib14}%
Arnold Walfisz,
{\sl Weylsche Exponentialsummen in der neueren Zahlentheorie}
(Berlin:
VEB Deutscher Verlag der Wissenschaften,
1963).
\endref

\beginref \parindent=20pt\item{\bref{15}}
\hldest{xyz}{}{bib15}%
Edmund Taylor Whittaker
and George Neville Watson,
{\sl A Course of Modern Analysis, 4th ed.}
(Cambridge:
Cambridge University Press,
1996).
\endref

\goodbreak\
        \bye